\documentclass{article}
\usepackage{amsfonts}
\usepackage{amsmath, amssymb}

\begin{document}

\begin{center}
{\large Some applications of MV algebras }

\begin{equation*}
\end{equation*}

Cristina FLAUT

\begin{equation*}
\end{equation*}
\end{center}

\textbf{Abstract. }{\small In this paper, some properties and applications
of MV-algebras are provided. We define a Fibonacci sequence in an MV-algebra
and we prove that such a stationary sequence gives us an idempotent element.
Taking into account of the representation of a finite MV-agebra, by using
Boolean elements of this algebra, we prove that such a sequence is always
stationary and is not periodic as in the situation when such a sequence is
studied on }t{\small he group }$\left( \mathbb{Z}_{n},+\right) ${\small ,
the group of integers modulo }$n$. {\small Moreover, as an application in
Coding Theory, to a Boolean algebra is attached a binary block code and is
proved that, under some conditions, the converse is also true}.

\begin{equation*}
\end{equation*}


\begin{equation*}
\end{equation*}

1. \textbf{Introduction}%
\begin{equation*}
\end{equation*}

C. C. Chang, in the paper [CHA; 58], introduced MV-algebras as a
generalization of Boolean algebras. In the last decades, numerous papers
have been devoted to the study of the properties and the applications of
MV-algebras, Wajsberg algebras and Boolean algebras. A recent applications
of these algebras were provided in [FV; 20], where were presented an
algorithm which can built Wajsberg algebras staring from a binary block
codes.

In this paper, some properties and applications of MV-algebras are given.
\qquad \qquad We define a Fibonacci sequence in such an algebra and we prove
that a Fibonacci stationary sequence gives us an idempotent element. Taking
into account of the representation of a finite MV-agebra, by using Boolean
elements of this algebra, we prove that such a sequence is always
stationary. This result is interesting, since in the group $\left( \mathbb{Z}%
_{n},+\right) $, the group of integers modulo $n$, the Fibonacci sequences
are periodic, with period given by Pisano period. From here, we can remark
that the representation of a finite MV-algebra can gives us a method for
finding new interesting results.

In Section 3, some examples of MV algebras and the number of their
idempotent elements are given. In Section 4, was provided an algorithm to
built a Boolean algebra of order $2^{k+1}$, starting from a Boolean algebra
of order two. Moreover, as an application in Coding Theory, to a Boolean
algebra is attached a binary block code and is proved that, under some
conditions, the converse is also true.\medskip 

\textbf{Definition 1.1.} The following ordered set $\left( \mathcal{L},\leq
\right) $ is called \textit{lattice} if and only for all elements$~x,y\in 
\mathcal{L}$ there are their supremum and infimum elements, $sup\{x,y\}$ and 
$inf\{x,y\}$, denoted by

\begin{equation*}
sup\{x,y\}=x\vee y\text{ and }inf\{x,y\}=x\wedge y.
\end{equation*}

The lattice $\left( \mathcal{L},\leq \right) $ is called a \textit{%
distributive lattice } if and only if, for each elements $x,y,z$, we have
the following relations: 
\begin{equation*}
x\vee \left( y\wedge z\right) =(x\vee y)\wedge (x\vee z),
\end{equation*}%
\begin{equation*}
x\wedge \left( y\vee z\right) =(x\wedge y)\vee (x\wedge z).
\end{equation*}

A lattice $\left( \mathcal{L},\leq \right) $ is called a \textit{bounded
lattice} if and only if there are an element $0$ being the least element in $%
\mathcal{L}$ and an element $1$ being the greatest element in $\mathcal{L}.$

In a lattice $\left( \mathcal{L},\leq \right) $ an element $x\in \mathcal{L}$
has a \textit{complement }if and only if there is an element $y\in \mathcal{L%
}$ satisfying the following relations: 
\begin{equation*}
x\vee y=1~\text{and}~x\wedge y=0.
\end{equation*}%
An element having a complement is called \textit{complemented}. We remark
that a complement of an element is not unique. If $\left( \mathcal{L},\leq
\right) $ is distributive, then each element has at most a complement.

A lattice $\left( \mathcal{L},\leq \right) $ is called a \textit{%
complemented lattice} if and only if it is a bounded lattice and each
element $x\in \mathcal{L}$ has a complement.\medskip

\textbf{Definition 1.2.} ([CHA; 58]) We consider an abelian monoid $\left(
X,\oplus ,0\right) $ equipped with an unary operation $\,"\lceil ",$ such
that the following conditions are satisfied:

i) $x\oplus \lceil 0=\lceil 0;$

ii) $\lceil (\lceil x)=x;$

iii) $\lceil \left( \lceil x\oplus y\right) \oplus y=$ $\lceil \left( \lceil
y\oplus x\right) \oplus x$, for all elements $x,y\in X$. This abelian monoid
is called an \textit{MV-algebra} and we denote it by $\left( X,\oplus
,\lceil ,0\right) .\medskip $

\textbf{Remark 1.3. }1) ([Mu; 07]) With the above notations, in an
MV-algebra we denote the constant element $\lceil 0$ with $1$, therefore

\begin{equation*}
1=\lceil 0.
\end{equation*}%
Considering the following multiplications

\begin{equation*}
x\odot y=\lceil \left( \lceil x\oplus \lceil y\right)
\end{equation*}%
and

\begin{equation*}
x\ominus y=x\odot \lceil y=\lceil \left( \lceil x\oplus y\right) ,
\end{equation*}%
we have that 
\begin{equation*}
x\oplus y=\lceil (\lceil x\odot \lceil y)\text{.}
\end{equation*}

2) ([COM; 00], Lemma 1.1.3) For each $x\in X$, the relations $x\oplus \lceil
x=1$ and $x\odot \lceil x=0$ are satisfied.$\medskip \medskip $

\textbf{Proposition 1.4. }([Mu; 07]) \textit{For the MV-algebra} $\left(
X,\oplus ,\lceil ,0\right) $ \textit{and }$x,y\in X$\textit{, the following
conditions are equivalent:}

\textit{a)~}$x\odot \lceil y=0;$

\textit{b)} $\lceil x\oplus y=\lceil 0=1;$

\textit{c)} $y=x\oplus \left( y\ominus x\right) =x\oplus \lceil \left(
\lceil y\oplus x\right) ;$

\textit{d)} \textit{An element} $w\in X$ \textit{such that}~$x\oplus w=y$ 
\textit{can be found}.\medskip $_{{}}\medskip $

\textbf{Definition 1.5.} ([Mu; 07]) We consider the MV-algebra $\left(
X,\oplus ,\lceil ,0\right) $. For $x,y\in X$, the following order relation
are defined on $X$:%
\begin{equation*}
x\leq y\text{ if and only if }\lceil x\oplus y=\lceil 0=1.
\end{equation*}

\textbf{Remark} \textbf{1.6.} i) From the above, we have that the definition
of the order relation on the MV-algebra $\left( X,\oplus ,\lceil ,0\right) $
can be done by using one of the equivalent conditions a)-d) from the above
proposition.

ii) ([COM; 00], Proposition 1.1.5) The order relation defined above gives us
a lattice structure on an MV-algebra:

$a)~x\vee y=(x\odot \lceil y)\oplus y=\left( x\ominus y\right) \oplus
y=\lceil \left( \lceil x\oplus y\right) \oplus y;$

$b)$ $x\wedge y=\lceil \left( \lceil x\vee \lceil y\right) =x\odot \left(
\lceil x\oplus y\right) $. We will denote this lattice with $\mathcal{L}%
\left( X\right) .\medskip $

\textbf{Definition 1.7. }([CHA; 58]) We consider an algebra $\left( W,\ast ,%
\overline{\phantom{x}},1\right) $ equipped with a binary operation $"\ast "$
and a unary operation $"\overline{\phantom{x}}"$ satisfying the following
conditions, for every $x,y,z\in W$:

i) $\left( x\ast y\right) \ast \left[ \left( y\ast z\right) \circ \left(
x\ast z\right) \right] =1;$

ii) $\left( x\ast y\right) \ast y=\left( y\ast x\right) \ast x;$

iii) $\left( \overline{x}\ast \overline{y}\right) \ast \left( y\ast x\right)
=1;$

iv) $1\ast x=x.$

This algebra is called a \textit{Wajsberg algebra.\medskip }

\textbf{Remark 1.8. }([COM; 00], Lemma 4.2.2 and Theorem 4.2.5)

a) For a Wajsberg algebra $\left( W,\ast ,\overline{\phantom{x}},1\right) $,
if we define the following multiplications 
\begin{equation*}
x\odot y=\overline{\left( x\ast \overline{y}\right) }
\end{equation*}%
and 
\begin{equation*}
x\oplus y=\overline{x}\ast y,
\end{equation*}%
for all $x,y\in W$, the obtained algebra $\left( W,\oplus ,\odot ,\overline{%
\phantom{x}},0,1\right) $ is an MV-algebra, with $0=\overline{1}$.

b) If on the MV-algebra $\left( X,\oplus ,\odot ,\lceil ,0,1\right) $ we
define the operation%
\begin{equation*}
x\ast y=\lceil x\oplus y,
\end{equation*}%
it results that $\left( X,\ast ,\lceil ,1\right) $ is a Wajsberg
algebra.\medskip

\textbf{Definition 1.9.} For the finite MV-algebras $\left( X,\ast ,%
\overline{\phantom{x}},0_{1}\right) $ and $\left( Y,\cdot ,^{\prime
},0_{2}\right) ,$ we define on their Cartesian product $Z=X\times Y$ the
following multiplication $"\Delta "$,%
\begin{equation}
\left( x_{1},y_{1}\right) \Delta \left( x_{2},y_{2}\right) =\left( x_{1}\ast
x_{2},y_{1}\cdot y_{2}\right) ,  \tag{1.1.}
\end{equation}%
The complement of the element $\left( x_{1},y_{1}\right) $ is $\rfloor
\left( x_{1},y_{1}\right) =$ $\left( \overline{x}_{1},x_{2}^{\prime }\right) 
$ and $0=\left( 0_{1},0_{2}\right) $. Therefore, by straightforward
calculation, we obtain that $\left( Z,\Delta ,\rfloor ,0\right) $ is also an
MV-algebra.\medskip

\textbf{Definition 1.10.}\ The algebra $\left( \mathcal{B},\vee \wedge
,\rceil ,0,1\right) $, equipped with two binary operations $\vee $ and $%
\wedge $ and a unary operation$\ \rceil $,$~$is called a \textit{Boolean
algebra} if and only if $\left( \mathcal{B},\vee \wedge \right) ~$is a
distributive and a complemented lattice with 
\begin{equation*}
b\vee \rceil b=1,
\end{equation*}%
\begin{equation*}
b\wedge \rceil b=0,
\end{equation*}%
for all elements $b\in \mathcal{B}$. The elements $0$ and $1$ are the least
and the greatest elements from the algebra $\mathcal{B}.\medskip $

\textbf{Remark 1.11.} Boolean algebras represent a particular case of
MV-algebras. Indeed, if $\left( \mathcal{B},\vee \wedge ,\rceil ,0,1\right) $
is a Boolean algebra, then can be easily checked that $\left( \mathcal{B}%
,\vee ,\rceil ,0\right) $ is an MV-algebra.\medskip

\textbf{Remark 1.12.} ([COM; 00], p.25 )

1) With the above notations, for each MV-algebra $\left( X,\oplus ,\lceil
,0\right) $, we have that $\mathcal{L}\left( X\right) $ is a distributive
lattice. For the algebra $X$ we will denote $\mathcal{B}\left( X\right) $ or 
$\mathcal{B}$($\mathcal{L}\left( X\right) )$ the set of all complemented
elements in $X$. The elements from $\mathcal{B}\left( X\right) $ are called
Boolean or idempotent elements.

2) Let $\left( X,\oplus ,\lceil ,0\right) $ be an MV-algebra. Therefore $%
x\in \mathcal{B}\left( X\right) $ if and only if $x\oplus y=x\vee y,$ for
all $y\in X$.\medskip

\textbf{Definition 1.13.} Let $(\mathcal{L},\vee \wedge )$ be a lattice with 
$0$ and $1$, the least and the greatest elements from $\mathcal{L}$. A
nonempty subset $\mathcal{I}\subseteq \mathcal{L}$ is called \textit{an} 
\textit{ideal} of the lattice $\mathcal{L}$ if and only if the following
conditions are satisfied:

a) $0\in \mathcal{I};$

b) If $x\in \mathcal{I}$ and $y\leq x,$ then $y\in \mathcal{I};$

c) If $x,y\in \mathcal{I}$, therefore $x\vee y\in \mathcal{I}.$

If $x\in \mathcal{L},$the set 
\begin{equation*}
\{z\in \mathcal{L}~/~z\leq x\}
\end{equation*}%
is called \textit{the principal ideal generated by} $x$ and will be denoted
by $(-\infty ,x].\medskip $

\textbf{Remark 1.14.}

i) ([COM; 00], Theorem 6.4.1) We consider $\left( X,\oplus ,\lceil ,0\right) 
$ an MV-algebra. For an element $\beta \in \mathcal{B}\left( X\right) $, we
have that $\left( (-\infty ,\beta ],\oplus ,\lceil ^{\beta },0\right) $ is
an MV-algebra, where $\lceil ^{\beta }x=\beta \wedge \lceil x$.

ii) ([COM; 00], Lemma 6.4.5) For the MV-algebra $\left( X,\oplus ,\lceil
,0\right) $ we consider the elements $x_{1},x_{2},...,x\,_{k}\in \mathcal{B}%
\left( X\right) -\{0,e\},k\geq 2$, such that

a) $x_{1}\vee x_{2}\vee ...\vee x_{k}=1;$

b) For $i\neq j$, we have $x_{i}\wedge x_{j}=0,i,j\in \{1,2,...,k\}.$

Therefore, we have that 
\begin{equation*}
X\simeq (-\infty ,x_{1}]\times (-\infty ,x_{2}]\times ...\times (-\infty
,x_{k}].
\end{equation*}%
iii\textbf{) }If $x_{1},x_{2},...,x\,_{k}\in \mathcal{B}\left( X\right)
-\{0,e\}$, the above decomposition is proper.

\begin{equation*}
\end{equation*}

\textbf{2}. \textbf{Main results}%
\begin{equation*}
\end{equation*}

In the last decades, a lot of papers have been devoted to the study of the
properties and applications of Fibonacci sequences in various algebraic
structure, as for example [HKN; 12], [KNS; 13], [Re; 13], etc. \ In the
following, we will prove that a Fibonacci sequence defined on a finite
MV-algebra is stationary and not periodic. This result is interesting
comparing with the behavior of such a sequence on the group $\left( \mathbb{Z%
}_{n},+\right) $, the group of integers modulo $n$, \ where the Fibonacci
sequences are periodic, with period given by Pisano period.

Let $X=\{x_{0}\leq x_{1}\leq ...\leq x_{n}\}$ be a finite totally ordered
set, with $x_{0}$ the minimum element and $x_{n}$ the maximum element. The
following multiplication $"\ast "$ are defined on $X$:%
\begin{equation}
\left\{ 
\begin{array}{c}
x_{i}\ast x_{j}=1\text{, if }x_{i}\leq x_{j}; \\ 
x_{i}\ast x_{j}=x_{n-i+j}\text{, otherwise;} \\ 
x_{0}=0,x_{n}=1,x\circ 0=\rfloor x.%
\end{array}%
\right.  \tag{2.1.}
\end{equation}%
It results that $\left( X,\ast ,\rfloor ,1\right) $ is a Wajsberg algebra.
As was remarked in [FRT; 84], Theorem 19, relation $\left( 2.1\right) $ give
us the only modality in which a Wajsberg algebra structure can be defined on
a finite totally ordered set, such that, on this algebra, the induced order
relation is given by $(2.1)$. Moreover, the relation $\rfloor x_{i}=x_{n-i}$
is fulfilled.\medskip

\textbf{Remark 2.1.} 1) Theorem 5.2, p. 43, from [HR; 99] tell us that an
MV-algebra is finite if and only if it is isomorphic to a finite product of
finite totally ordered MV-algebras. Using connections between MV-algebras
and Wajsberg algebras, it results that if \ $M=\left( X,\oplus ,\odot
,\lceil ,0,1\right) $ is a totally ordered MV-algebra, then the obtained
Wajsberg algebra, $\,W=\left( X,\ast ,\lceil ,1\right) $, is also totally
ordered. The converse of this statement is also true, since $x\ast y=\lceil
x\oplus y$ implies that $x\leq _{M}y$ if and only $x\leq _{W}y$.

If the number of elements is a finite MV-algebra or in a finite Wajsberg
algebra is a prime number, therefore these algebras are totally ordered
algebras. Remark 1.14, ii) gives us a similar result and, additionally, we
obtain that the sets $(-\infty ,x_{i}],i\in \{1,2,...,k\},$ from that
decomposition are totally ordered.\medskip

\textbf{Example 2.2. (}[WDH; 17], Example 3.3\textbf{) }We consider the
following MV-algebra $\left( X,\oplus ,\lceil ,0\right) $, with the
multiplication $\oplus $ and the operation $"\lceil "$ given in the below
tables:

\begin{center}
\begin{tabular}{l|llllll}
$\oplus $ & $0$ & $\alpha $ & $\beta $ & $\gamma $ & $\delta $ & $%
\varepsilon $ \\ \hline
$0$ & $0$ & $\alpha $ & $\beta $ & $\gamma $ & $\delta $ & $\varepsilon $ \\ 
$\alpha $ & $\alpha $ & $\gamma $ & $\delta $ & $\gamma $ & $\varepsilon $ & 
$\varepsilon $ \\ 
$\beta $ & $\beta $ & $\delta $ & $\beta $ & $\varepsilon $ & $\delta $ & $%
\varepsilon $ \\ 
$\gamma $ & $\gamma $ & $\gamma $ & $\varepsilon $ & $\gamma $ & $%
\varepsilon $ & $\varepsilon $ \\ 
$\delta $ & $\delta $ & $\varepsilon $ & $\delta $ & $\varepsilon $ & $%
\varepsilon $ & $\varepsilon $ \\ 
$\varepsilon $ & $\varepsilon $ & $\varepsilon $ & $\varepsilon $ & $%
\varepsilon $ & $\varepsilon $ & $\varepsilon $%
\end{tabular}
~\ 
\begin{tabular}{l|llllll}
$\lceil $ & $\theta $ & $\alpha $ & $\beta $ & $\gamma $ & $\delta $ & $%
\varepsilon $ \\ \hline
& $\varepsilon $ & $\delta $ & $\gamma $ & $\beta $ & $\alpha $ & $0$%
\end{tabular}%
.
\end{center}

We remark that $\beta \oplus \beta =\beta $ and $\gamma \oplus \gamma
=\gamma $. We have\newline
$\beta \vee \gamma =\lceil \left( \gamma \oplus \gamma \right) \oplus \gamma
=\beta \oplus \gamma =\varepsilon $ and \newline
$\beta \wedge \gamma =\lceil \left( \lceil \beta \vee \lceil \gamma \right)
=\lceil \left( \gamma \vee \beta \right) =\lceil \varepsilon =0$. Now, we
compute $(-\infty ,\beta ]$ and $(-\infty ,\gamma ]$. By using Remark 1.14,
ii), it results: $(-\infty ,\beta ]=\{0,\beta \},(-\infty ,\gamma
]=\{0,\alpha ,\gamma \}$. Therefore, $X\simeq (-\infty ,\beta ]\times
(-\infty ,\gamma ]$.\medskip

\textbf{Definition 2.3.} We consider $\left( X,\oplus ,\lceil ,0\right) $ an
MV-algebra. For $x,y\in X$, the following sequence 
\begin{equation*}
<x,y>=\{x,y,x\oplus y,y\oplus \left( x\oplus y\right)
,...,.u_{n},u_{n+1},u_{n+2},.....\}
\end{equation*}%
are defined, where $u_{0}=x,u_{1}=y$ and $u_{n+2}=u_{n}\oplus u_{n+1}$, for $%
k\in \mathbb{N}$. This sequence are called the \textit{Fibonacci sequence
attached to the elements} $x,y$. If there is a number $k\in \mathbb{N}$ such
that $u_{n}=u_{n+1}=u_{n+2}=....$, for all $n\geq k$, then the sequence $%
<x,y>$ is called $k$\textit{-stationary}.\medskip

\textbf{Proposition} \textbf{2.4}. \textit{For the MV-algebra } $\left(
X,\oplus ,\lceil ,0\right) $ \textit{we consider} $x,y\in X$. \textit{If the
sequence} $<x,y>$ \textit{is} $k$\textit{-stationary therefore} $u_{k}$ 
\textit{is idempotent.\medskip }

\textbf{Proof.} Let $k\in \mathbb{N}$ be a natural number such that the
sequence $<x,y>$ is $k$-stationary. Therefore, $u_{k}=u_{k+1}=u_{k+2}=...$.
We have $u_{2}=x\oplus y,u_{3}=y\oplus \left( x\oplus y\right) =x\oplus
2y,2y=y\oplus y,u_{4}=2x\oplus 3y$, etc. It results that 
\begin{equation*}
u_{n}=f_{n-1}x\oplus f_{n}y,
\end{equation*}%
where $(f_{n})_{n\in \mathbb{N}}$ is the Fibonacci sequence 
\begin{equation*}
f_{0}=0,f_{1}=1,f_{n+1}=f_{n}+f_{n-1},n\in \mathbb{N}\text{, }n\geq 1.
\end{equation*}

Since $u=u_{k}=u_{k+1}=u_{k+2}=...$, we have\newline
$u=f_{k-1}x\oplus f_{k}y=f_{k}x\oplus f_{k+1}y=f_{k+1}x\oplus f_{k+2}y$. 
\newline
Therefore, $u=u_{k+2}=f_{k+1}x\oplus f_{k+2}y=(f_{k-1}x\oplus f_{k}y)\oplus
(f_{k}x\oplus f_{k+1}y)=u\oplus u$. It results that $u=u_{k}$ is an
idempotent (Boolean) element.$_{{}}\medskip $

\textbf{Proposition 2.5.} \textit{Let} $\left( X,\oplus ,\lceil ,0\right) $ 
\textit{be an MV-algebra and} $x,y\in X$. \textit{If the sequence} $<x,y>$ 
\textit{is} $2$\textit{-stationary for all }$x,y\in X$\textit{, therefore} $%
X $ \textit{is a Boolean algebra.\medskip }

\textbf{Proof.} \ We have that $x\oplus y$ is idempotent, therefore $x\oplus
y=(x\oplus y)\oplus (x\oplus y)$, for all $x,y\in X$. If we take $y=0,$
therefore $x=x\oplus x$, for all $x\in X$. It results that $X$ is a Boolean
algebra.$_{{}}\medskip $

\textbf{Proposition 2.6.} \textit{Let} $\left( X,\oplus ,\lceil ,0\right) $ 
\textit{be a finite MV-algebra and} $x,y\in X$. \textit{Therefore the
sequence} $<x,y>$ \textit{is} \textit{stationary for all }$x,y\in X$\textit{%
.\medskip }

\textbf{Proof.} Since a finite MV-algebra is isomorphic to a finite product
of finite totally ordered MV-algebras, it is enough to proof this result in
the case of the finite totally ordered algebras. Let $X=%
\{x_{0},x_{1},...,x_{n}\}$ be a finite totally ordered MV-algebra. From
relation $\left( 2.1\right) $ and since $\lceil x_{i}=x_{n-i}$, we have that
the multiplication $"\oplus "$ is given by the following formulae:

\begin{equation*}
\left\{ 
\begin{array}{c}
x_{i}\oplus x_{j}=1\text{, if }i+j>n; \\ 
x_{i}\oplus x_{j}=x_{i+j}\text{, if }i+j\leq n\text{;} \\ 
x_{0}=0,x_{n}=1,x\oplus 0=x.%
\end{array}%
\right. .
\end{equation*}%
\textbf{Case 1.} Let $x,y\in X,x\neq 0,y\neq 0,$ and the sequence\newline
\begin{equation*}
<x,y>=\{x,y,x\oplus y,y\oplus \left( x\oplus y\right)
,...,.u_{n},u_{n+1},u_{n+2},.....\},
\end{equation*}%
where $u_{0}=x,u_{1}=y$ and $u_{n+2}=u_{n}\oplus u_{n+1}$, for $k\in \mathbb{%
N}$. It is clear that $u_{2}=x\oplus y>x$ and $u_{2}=x\oplus y>y$. If $%
x\oplus y=x$ or $x\oplus y=y,$ therefore $y=\theta $ or $x=\theta $, false.
Assuming $u_{2}=x\oplus y>x$, we have that\newline
$u_{3}=y\oplus \left( x\oplus y\right) =\left( x\oplus y\right) \oplus
y>x\oplus y=u_{2},$\newline
$u_{4}=(x\oplus y)\oplus \left( y\oplus \left( x\oplus y\right) \right)
>y\oplus \left( x\oplus y\right) =u_{3},$\newline
$u_{5}=\left( x\oplus y\right) \oplus \left( (x\oplus y)\oplus \left(
y\oplus \left( x\oplus y\right) \right) \right) >\left( (x\oplus y)\oplus
\left( y\oplus \left( x\oplus y\right) \right) \right) =u_{4}$, etc. It
results that the obtained increased sequence is stationary, since the set $X$
is finite. Therefore, there is $k\in \mathbb{N}$ such that $%
u_{k}=1=u_{k+1}=u_{k+2}=...$. We get

\begin{equation*}
<x,y>=\{x,y,x\oplus y,y\oplus \left( x\oplus y\right)
,...,.u_{k-1},1,1,1,1,.....\}.
\end{equation*}

\textbf{Case 2.} Let $x,y\in X,y=0$. We obtain the sequence $<x,y>$ with $%
u_{0}=x,u_{1}=0,u_{2}=x,u_{3}=x,u_{4}=$ $x\oplus x>x$. We apply the Case 1,
obtaining a stationary sequence, that means a number $k\in \mathbb{N}$ such
that $u_{k}=1=u_{k+1}=u_{k+2}=...$.

Now, by using Definition 1.9, it is clear that in a finite MV-algebra\textit{%
\ }the sequence $<x,y>$ is stationary for all\textit{\ }$x,y\in X$\textit{. }%
$_{{}}\medskip $

\textbf{Remark 2.7.}

1) The above result is not true for infinite MV-algebras, as can be easily
seen by using the famous Chang's MV-algebra.

2) Let $\left( X,\oplus ,\lceil ,0\right) $ be a finite MV-algebra and $%
x,y\in X$. Since from the above it results that the algebra $X$ is $k$%
-stationary, we have that $u_{n}=a,$ for all $n\geq k$. We consider the
following map 
\begin{equation*}
\lambda :X\times X\rightarrow X,\lambda \left( x,y\right) =u_{k}.
\end{equation*}%
In MV-algebra $\left( X,\oplus ,\lceil ,0\right) $ given in the Example 2.2,
we have that $\lambda \left( \gamma ,0\right) =\gamma $, $\lambda \left(
\gamma ,\delta \right) =\varepsilon $, $\lambda \left( \alpha ,\beta \right)
=\varepsilon $. Indeed, \newline
-the sequence $\left[ \gamma ,0\right] =\gamma ,0,\gamma ,\gamma ,....$ is $2
$-stationary;\newline
-the sequence $\left[ \gamma ,\delta \right] =\gamma ,\delta ,\varepsilon
,\varepsilon ,...$ is $2$-stationary;\newline
-the sequence $\left[ \alpha ,\beta \right] =\alpha ,\beta ,\delta ,\delta
,\varepsilon ,\varepsilon ,...$ is $4$-stationary.\medskip 

\textbf{Remark 2.8.} If $\left( X,\oplus ,\lceil ,0\right) $ is a finite
MV-algebra such that $X\simeq (-\infty ,x_{1}]\times (-\infty ,x_{2}]\times
...\times (-\infty ,x_{k}]$ and the sets $(-\infty ,x_{i}]=\{0,x_{i}\}$,
that means it has two elements for all $i\in \{1,2,...,k\}$, it results that 
$X$ is a Boolean algebra. Indeed, using above results, we have that all
elements in $X$ have the form $\left( \alpha _{1},\alpha _{2},...,\alpha
_{k}\right) $, where $\alpha _{i}$ $\in $ $(-\infty ,x_{i}]$. From here, we
get the known result that a finite Boolean algebra has $2^{k}$ elements.

\begin{equation*}
\end{equation*}

\textbf{3.} \textbf{Examples}%
\begin{equation*}
\end{equation*}

In papers [BV; 10], [FHSV; 20] was presented classification of MV-algebras
by using different algorithms. In [FV; 20] was presented an application of
these algebras in Coding Theory. In the following, by using Remark 1.14,
ii), which give an alternative method to characterize MV-algebras and, as a
consequence, Boolean algebras, we will give some examples of MV-algebras and
Boolean algebras.\medskip

\textbf{Example 3.1.} Using examples from [FHSV; 20], Section 4.1, we
consider $W=\{0\leq \alpha \leq \beta \leq \varepsilon \},$ a totally
ordered set. We define on $W$ two multiplications $\Delta _{0}^{4}$ and $%
\oplus _{0}^{4}~$given in the below tables. With multiplication $\Delta
_{0}^{4}$ with $\overline{\alpha }=\beta $ and $\overline{\beta }=\alpha $, $%
W$ becomes a Wajsberg algebra. The associated MV-algebra is obtained with
multiplication $\oplus _{0}^{4}$. It results that in the MV-algebra
structure, the only idempotent elements are $0$ and $\varepsilon $.

\begin{equation*}
\begin{tabular}{l|llll}
$\Delta _{0}^{4}$ & $0$ & $\alpha $ & $\beta $ & $\varepsilon $ \\ \hline
$0$ & $\varepsilon $ & $\varepsilon $ & $\varepsilon $ & $\varepsilon $ \\ 
$\alpha $ & $\beta $ & $\varepsilon $ & $\varepsilon $ & $\varepsilon $ \\ 
$\beta $ & $\alpha $ & $\beta $ & $\varepsilon $ & $\varepsilon $ \\ 
$\varepsilon $ & $0$ & $\alpha $ & $\beta $ & $\varepsilon $%
\end{tabular}%
\ \ \ \ \ \ \ 
\begin{tabular}{l|llll}
$\oplus _{0}^{4}$ & $0$ & $\alpha $ & $\beta $ & $\varepsilon $ \\ \hline
$0$ & $0$ & $\alpha $ & $\beta $ & $\varepsilon $ \\ 
$\alpha $ & $\alpha $ & $\beta $ & $\varepsilon $ & $\varepsilon $ \\ 
$\beta $ & $\beta $ & $\varepsilon $ & $\varepsilon $ & $\varepsilon $ \\ 
$\varepsilon $ & $\varepsilon $ & $\varepsilon $ & $\varepsilon $ & $%
\varepsilon $%
\end{tabular}%
.
\end{equation*}

We consider now the set $W=\{0,\alpha ,\beta ,\varepsilon \},$\ partially
ordered. On $W$ we define two multiplications $\Delta _{11}^{4}$ and $\oplus
_{11}^{4}$given in the below tables. With multiplication $\Delta _{11}^{4},$ 
$W$ becomes a Wajsberg algebra. The associated MV-algebra is obtained with
multiplication $\oplus _{11}^{4}$.

\begin{equation}
\begin{tabular}{l|llll}
$\Delta _{11}^{4}$ & $0$ & $\alpha $ & $\beta $ & $\varepsilon $ \\ \hline
$0$ & $\varepsilon $ & $\varepsilon $ & $\varepsilon $ & $\varepsilon $ \\ 
$\alpha $ & $\beta $ & $\varepsilon $ & $\beta $ & $\varepsilon $ \\ 
$\beta $ & $\alpha $ & $\alpha $ & $\varepsilon $ & $\varepsilon $ \\ 
$\varepsilon $ & $0$ & $\alpha $ & $\beta $ & $\varepsilon $%
\end{tabular}%
\ ~~\ \ 
\begin{tabular}{l|llll}
$\oplus _{11}^{4}$ & $0$ & $\alpha $ & $\beta $ & $\varepsilon $ \\ \hline
\multicolumn{1}{l|}{$0$} & $0$ & $\alpha $ & $\beta $ & $\varepsilon $ \\ 
\multicolumn{1}{l|}{$\alpha $} & $\alpha $ & $\alpha $ & $\varepsilon $ & $%
\varepsilon $ \\ 
\multicolumn{1}{l|}{$\beta $} & $\beta $ & $\varepsilon $ & $\beta $ & $%
\varepsilon $ \\ 
\multicolumn{1}{l|}{$\varepsilon $} & $\varepsilon $ & $\varepsilon $ & $%
\varepsilon $ & $\varepsilon $%
\end{tabular}%
\text{.}  \tag{3.1.}
\end{equation}

In the above obtained MV-algebra structure all elements are idempotent. For $%
\alpha $ and $\beta $ we remark that $\alpha \vee \beta =\varepsilon $ and $%
\alpha \wedge \beta =0$. Therefore, as MV-algebra, $W$ $\simeq (-\infty
,\alpha ]\times (-\infty ,\beta ]$, where $(-\infty ,\alpha ]=\{0,\alpha \}$
and $(-\infty ,\beta ]=\{0,\beta \}$. From here, we obtain that there exist
only two non-isomorphic MV-algebras of order $4$. Thus, we can remark that
in an MV-algebra of order $4$ we can have only $0$ or $2$ proper
idempotents. The MV-algebra $\left( W,\oplus _{11}^{4}\right) $ is the only
Boolean algebra of order $4$. We denote this algebra with $\mathcal{B}_{4}$%
.\medskip

\textbf{Example 3.2}. Using examples from [FHSV; 20], Section 4.2, we
consider $W=\{0\leq \alpha \leq \beta \leq \gamma \leq \delta \leq
\varepsilon \},$ a totally ordered set. We define on $W$ two multiplications 
$\Delta _{0}^{6}$ and $\oplus _{0}^{6}$given in the below tables. With
multiplication $\Delta _{0}^{6},$ with $\overline{\alpha }=\delta $, $%
\overline{\beta }=\gamma $, $\overline{\gamma }=\beta $, $\overline{\delta }%
=\alpha $, $W$ becomes a Wajsberg algebra. The associated MV-algebra is
obtained with multiplication $\oplus _{0}^{6}$. We can see that in the
MV-algebra structure, the only idempotent elements are $0$ and $\varepsilon $%
.

\begin{equation*}
\begin{tabular}{l|llllll}
$\Delta _{0}^{6}$ & $0$ & $\alpha $ & $\beta $ & $\gamma $ & $\delta $ & $%
\varepsilon $ \\ \hline
$0$ & $\varepsilon $ & $\varepsilon $ & $\varepsilon $ & $\varepsilon $ & $%
\varepsilon $ & $\varepsilon $ \\ 
$\alpha $ & $\delta $ & $\varepsilon $ & $\varepsilon $ & $\varepsilon $ & $%
\varepsilon $ & $\varepsilon $ \\ 
$\beta $ & $\gamma $ & $\delta $ & $\varepsilon $ & $\varepsilon $ & $%
\varepsilon $ & $\varepsilon $ \\ 
$\gamma $ & $\beta $ & $\gamma $ & $\delta $ & $\varepsilon $ & $\varepsilon 
$ & $\varepsilon $ \\ 
$\delta $ & $\alpha $ & $\beta $ & $\gamma $ & $\delta $ & $\varepsilon $ & $%
\varepsilon $ \\ 
$\varepsilon $ & $0$ & $\alpha $ & $\beta $ & $\gamma $ & $\delta $ & $%
\varepsilon $%
\end{tabular}%
\text{ \ \ }%
\begin{tabular}{l|llllll}
$\oplus _{0}^{6}$ & $0$ & $\alpha $ & $\beta $ & $\gamma $ & $\delta $ & $%
\varepsilon $ \\ \hline
$0$ & $0$ & $\alpha $ & $\beta $ & $\gamma $ & $\delta $ & $\varepsilon $ \\ 
$\alpha $ & $\alpha $ & $\beta $ & $\gamma $ & $\delta $ & $\varepsilon $ & $%
\varepsilon $ \\ 
$\beta $ & $\beta $ & $\gamma $ & $\delta $ & $\varepsilon $ & $\varepsilon $
& $\varepsilon $ \\ 
$\gamma $ & $\gamma $ & $\delta $ & $\varepsilon $ & $\varepsilon $ & $%
\varepsilon $ & $\varepsilon $ \\ 
$\delta $ & $\delta $ & $\varepsilon $ & $\varepsilon $ & $\varepsilon $ & $%
\varepsilon $ & $\varepsilon $ \\ 
$\varepsilon $ & $\varepsilon $ & $\varepsilon $ & $\varepsilon $ & $%
\varepsilon $ & $\varepsilon $ & $\varepsilon $%
\end{tabular}%
.
\end{equation*}

We consider now the set $W=\{0,\alpha ,\beta ,\gamma ,\delta ,\varepsilon
\}, $partially ordered. On $W$ we define two multiplications $\Delta
_{11}^{6}$ and $\oplus _{11}^{6}$given in the below tables. With
multiplication $\Delta _{11}^{6},$ $W$ becomes a Wajsberg algebra. The
associated MV-algebra is obtained with the multiplication $\oplus _{11}^{6}$.

\begin{equation*}
\begin{tabular}{l|llllll}
$\Delta _{11}^{6}$ & $0$ & $\alpha $ & $\beta $ & $\gamma $ & $\delta $ & $%
\varepsilon $ \\ \hline
$0$ & $\varepsilon $ & $\varepsilon $ & $\varepsilon $ & $\varepsilon $ & $%
\varepsilon $ & $\varepsilon $ \\ 
$\alpha $ & $\delta $ & $\varepsilon $ & $\varepsilon $ & $\delta $ & $%
\varepsilon $ & $\varepsilon $ \\ 
$\beta $ & $\gamma $ & $\delta $ & $\varepsilon $ & $\gamma $ & $\delta $ & $%
\varepsilon $ \\ 
$\gamma $ & $\beta $ & $\beta $ & $\beta $ & $\varepsilon $ & $\varepsilon $
& $\varepsilon $ \\ 
$\delta $ & $\alpha $ & $\beta $ & $\beta $ & $\delta $ & $\varepsilon $ & $%
\varepsilon $ \\ 
$\varepsilon $ & $0$ & $\alpha $ & $\beta $ & $\gamma $ & $\delta $ & $%
\varepsilon $%
\end{tabular}%
\text{ \ \ }%
\begin{tabular}{lllllll}
$\oplus _{11}^{6}$ & \multicolumn{1}{|l}{$0$} & $\alpha $ & $\beta $ & $%
\gamma $ & $\delta $ & $\varepsilon $ \\ \hline
\multicolumn{1}{l|}{$0$} & \multicolumn{1}{l|}{$0$} & $\alpha $ & $\beta $ & 
$\gamma $ & $\delta $ & $\varepsilon $ \\ 
\multicolumn{1}{l|}{$\alpha $} & $\alpha $ & $\beta $ & $\beta $ & $\delta $
& $\varepsilon $ & $\varepsilon $ \\ 
\multicolumn{1}{l|}{$\beta $} & $\beta $ & $\beta $ & $\beta $ & $%
\varepsilon $ & $\varepsilon $ & $\varepsilon $ \\ 
\multicolumn{1}{l|}{$\gamma $} & $\gamma $ & $\delta $ & $\varepsilon $ & $%
\gamma $ & $\delta $ & $\varepsilon $ \\ 
\multicolumn{1}{l|}{$\delta $} & $\delta $ & $\varepsilon $ & $\varepsilon $
& $\delta $ & $\varepsilon $ & $\varepsilon $ \\ 
\multicolumn{1}{l|}{$\varepsilon $} & \multicolumn{1}{l|}{$\varepsilon $} & $%
\varepsilon $ & $\varepsilon $ & $\varepsilon $ & $\varepsilon $ & $%
\varepsilon $%
\end{tabular}%
.
\end{equation*}

We remark that in the MV-algebra structure the idempotent elements are $%
\{0,\beta ,\gamma ,\varepsilon \}$. For $\beta $ and $\gamma $ we remark
that $\beta \vee \gamma =\varepsilon $ and $\beta \wedge \gamma =0$.
Therefore, as MV-algebra, $W$ $\simeq (-\infty ,\beta ]\times (-\infty
,\gamma ]$, where $(-\infty ,\beta ]=\{0,\alpha ,\beta \}$, $(-\infty
,\gamma ]=\{0,\gamma \}$. From here, we get that there are only two
non-isomorphic MV-algebras of order $6$. Thus, we can remark that in an
MV-algebra of order $6$ we can have only $0$ or $2$ proper
idempotents.\medskip

\textbf{Example 3.3.} Using examples from [FHSV; 20], Section 4.3, we
consider \newline
$W=\{0\leq \alpha \leq \beta \leq \gamma \leq \tau \leq \upsilon \leq \rho
\leq \varepsilon \}$ a totally ordered set. On $W$ we define two
multiplications $\Delta _{0}^{8}$ and $\oplus _{0}^{8}$, given in the below
tables. With multiplication $\Delta _{0}^{8},$ with $\overline{\alpha }=\rho 
$, $\overline{\beta }=\upsilon $, $\overline{\gamma }=\tau $, $W$ becomes a
Wajsberg algebra. The associated MV-algebra is obtained with multiplication $%
\oplus _{0}^{8}$. We remark that in the MV-algebra structure, the only
idempotent elements are $0$ and $\varepsilon $.

\begin{equation*}
\begin{tabular}{l|llllllll}
$\Delta _{0}^{8}$ & $0$ & $\alpha $ & $\beta $ & $\gamma $ & $\tau $ & $%
\upsilon $ & $\rho $ & $\varepsilon $ \\ \hline
$0$ & $\varepsilon $ & $\varepsilon $ & $\varepsilon $ & $\varepsilon $ & $%
\varepsilon $ & $\varepsilon $ & $\varepsilon $ & $\varepsilon $ \\ 
$\alpha $ & $\rho $ & $\varepsilon $ & $\varepsilon $ & $\varepsilon $ & $%
\varepsilon $ & $\varepsilon $ & $\varepsilon $ & $\varepsilon $ \\ 
$\beta $ & $\upsilon $ & $\rho $ & $\varepsilon $ & $\varepsilon $ & $%
\varepsilon $ & $\varepsilon $ & $\varepsilon $ & $\varepsilon $ \\ 
$\gamma $ & $\tau $ & $\upsilon $ & $\rho $ & $\varepsilon $ & $\varepsilon $
& $\varepsilon $ & $\varepsilon $ & $\varepsilon $ \\ 
$\tau $ & $\gamma $ & $\tau $ & $\upsilon $ & $\rho $ & $\varepsilon $ & $%
\varepsilon $ & $\varepsilon $ & $\varepsilon $ \\ 
$\upsilon $ & $\beta $ & $\gamma $ & $\tau $ & $\upsilon $ & $\rho $ & $%
\varepsilon $ & $\varepsilon $ & $\varepsilon $ \\ 
$\rho $ & $\alpha $ & $\beta $ & $\gamma $ & $\tau $ & $\upsilon $ & $\rho $
& $\varepsilon $ & $\varepsilon $ \\ 
$\varepsilon $ & $0$ & $\alpha $ & $\beta $ & $\gamma $ & $\tau $ & $%
\upsilon $ & $\rho $ & $\varepsilon $%
\end{tabular}%
~\ 
\begin{tabular}{l|llllllll}
$\oplus _{0}^{8}$ & $0$ & $\alpha $ & $\beta $ & $\gamma $ & $\tau $ & $%
\upsilon $ & $\rho $ & $\varepsilon $ \\ \hline
\multicolumn{1}{l|}{$0$} & $0$ & $\alpha $ & $\beta $ & $\gamma $ & $\tau $
& $\upsilon $ & $\rho $ & $\varepsilon $ \\ 
\multicolumn{1}{l|}{$\alpha $} & $\alpha $ & $\beta $ & $\gamma $ & $\tau $
& $\upsilon $ & $\rho $ & $\varepsilon $ & $\varepsilon $ \\ 
\multicolumn{1}{l|}{$\beta $} & $\beta $ & $\gamma $ & $\tau $ & $\upsilon $
& $\rho $ & $\varepsilon $ & $\varepsilon $ & $\varepsilon $ \\ 
\multicolumn{1}{l|}{$\gamma $} & $\gamma $ & $\tau $ & $\upsilon $ & $\rho $
& $\varepsilon $ & $\varepsilon $ & $\varepsilon $ & $\varepsilon $ \\ 
\multicolumn{1}{l|}{$\tau $} & $\tau $ & $\upsilon $ & $\rho $ & $%
\varepsilon $ & $\varepsilon $ & $\varepsilon $ & $\varepsilon $ & $%
\varepsilon $ \\ 
\multicolumn{1}{l|}{$\upsilon $} & $\upsilon $ & $\rho $ & $\varepsilon $ & $%
\varepsilon $ & $\varepsilon $ & $\varepsilon $ & $\varepsilon $ & $%
\varepsilon $ \\ 
\multicolumn{1}{l|}{$\rho $} & $\rho $ & $\varepsilon $ & $\varepsilon $ & $%
\varepsilon $ & $\varepsilon $ & $\varepsilon $ & $\varepsilon $ & $%
\varepsilon $ \\ 
\multicolumn{1}{l|}{$\varepsilon $} & $\varepsilon $ & $\varepsilon $ & $%
\varepsilon $ & $\varepsilon $ & $\varepsilon $ & $\varepsilon $ & $%
\varepsilon $ & $\varepsilon $%
\end{tabular}%
\end{equation*}

We consider now the set $W=\{0,\alpha ,\beta ,\gamma ,\tau ,\upsilon ,\rho
,\varepsilon \},$ partially ordered. On $W$ we define two multiplications $%
\Delta _{11}^{8}$ and $\oplus _{11}^{8}$, given in the below tables. With
multiplication $\Delta _{11}^{8},$ $W$ becomes a Wajsberg algebra. The
associated MV-algebra is obtained with multiplication $\oplus _{11}^{8}.$

\begin{equation*}
\begin{tabular}{l|llllllll}
$\Delta _{11}^{8}$ & $0$ & $\alpha $ & $\beta $ & $\gamma $ & $\tau $ & $%
\upsilon $ & $\rho $ & $\varepsilon $ \\ \hline
$0$ & $\varepsilon $ & $\varepsilon $ & $\varepsilon $ & $\varepsilon $ & $%
\varepsilon $ & $\varepsilon $ & $\varepsilon $ & $\varepsilon $ \\ 
$\alpha $ & $\rho $ & $\varepsilon $ & $\rho $ & $\varepsilon $ & $\rho $ & $%
\varepsilon $ & $\rho $ & $\varepsilon $ \\ 
$\beta $ & $\upsilon $ & $\upsilon $ & $\varepsilon $ & $\varepsilon $ & $%
\varepsilon $ & $\varepsilon $ & $\varepsilon $ & $\varepsilon $ \\ 
$\gamma $ & $\tau $ & $\upsilon $ & $\rho $ & $\varepsilon $ & $\rho $ & $%
\varepsilon $ & $\rho $ & $\varepsilon $ \\ 
$\tau $ & $\gamma $ & $\gamma $ & $\upsilon $ & $\upsilon $ & $\varepsilon $
& $\varepsilon $ & $\varepsilon $ & $\varepsilon $ \\ 
$\upsilon $ & $\beta $ & $\gamma $ & $\tau $ & $\upsilon $ & $\tau $ & $%
\varepsilon $ & $\rho $ & $\varepsilon $ \\ 
$\rho $ & $\alpha $ & $\alpha $ & $\gamma $ & $\gamma $ & $\upsilon $ & $%
\upsilon $ & $\varepsilon $ & $\varepsilon $ \\ 
$\varepsilon $ & $0$ & $\alpha $ & $\beta $ & $\gamma $ & $\tau $ & $%
\upsilon $ & $\rho $ & $\varepsilon $%
\end{tabular}%
\text{ \ }%
\begin{tabular}{l|llllllll}
$\oplus _{11}^{8}$ & $0$ & $\alpha $ & $\beta $ & $\gamma $ & $\tau $ & $%
\upsilon $ & $\rho $ & $\varepsilon $ \\ \hline
\multicolumn{1}{l|}{$0$} & $0$ & $\alpha $ & $\beta $ & $\gamma $ & $\tau $
& $\upsilon $ & $\rho $ & $\varepsilon $ \\ 
\multicolumn{1}{l|}{$\alpha $} & $\alpha $ & $\alpha $ & $\gamma $ & $\gamma 
$ & $\upsilon $ & $\upsilon $ & $\varepsilon $ & $\varepsilon $ \\ 
\multicolumn{1}{l|}{$\beta $} & $\beta $ & $\gamma $ & $\tau $ & $\upsilon $
& $\tau $ & $\varepsilon $ & $\rho $ & $\varepsilon $ \\ 
\multicolumn{1}{l|}{$\gamma $} & $\gamma $ & $\gamma $ & $\upsilon $ & $%
\upsilon $ & $\varepsilon $ & $\varepsilon $ & $\varepsilon $ & $\varepsilon 
$ \\ 
\multicolumn{1}{l|}{$\tau $} & $\tau $ & $\upsilon $ & $\rho $ & $%
\varepsilon $ & $\rho $ & $\varepsilon $ & $\rho $ & $\varepsilon $ \\ 
\multicolumn{1}{l|}{$\upsilon $} & $\upsilon $ & $\upsilon $ & $\varepsilon $
& $\varepsilon $ & $\varepsilon $ & $\varepsilon $ & $\varepsilon $ & $%
\varepsilon $ \\ 
\multicolumn{1}{l|}{$\rho $} & $\rho $ & $\varepsilon $ & $\rho $ & $%
\varepsilon $ & $\rho $ & $\varepsilon $ & $\rho $ & $\varepsilon $ \\ 
\multicolumn{1}{l|}{$\varepsilon $} & $\varepsilon $ & $\varepsilon $ & $%
\varepsilon $ & $\varepsilon $ & $\varepsilon $ & $\varepsilon $ & $%
\varepsilon $ & $\varepsilon $%
\end{tabular}%
\end{equation*}

In this case, we remark that in this MV-algebra structure the idempotent
elements are $\{0,\alpha ,\rho ,\varepsilon \}$. For $\alpha $ and $\rho $
we remark that $\alpha \vee \rho =\varepsilon $ and $X\wedge \rho =0$.
Therefore, as MV-algebra, $W$ $\simeq (-\infty ,\alpha ]\times (-\infty
,\rho ]$, where $(-\infty ,\alpha ]=\{0,\alpha \}$, $(-\infty ,\rho
]=\{0,\beta ,\tau ,\rho \}$.

In the following, we consider the set $W=\{0,\alpha ,\beta ,\gamma ,\tau
,\upsilon ,\rho ,\varepsilon \}$,$~$partially ordered, where we define two
multiplications $\Delta _{21}^{8}$ and $\oplus _{21}^{8}$, given in the
below tables. With multiplication $\Delta _{21}^{8},$ $W$ becomes a Wajsberg
algebra. The associated MV-algebra is obtained with multiplication $\oplus
_{21}^{8}.$

\begin{equation}
\begin{tabular}{l|llllllll}
$\Delta _{21}^{8}$ & $0$ & $\alpha $ & $\beta $ & $\gamma $ & $\tau $ & $%
\upsilon $ & $\rho $ & $\varepsilon $ \\ \hline
$0$ & $\varepsilon $ & $\varepsilon $ & $\varepsilon $ & $\varepsilon $ & $%
\varepsilon $ & $\varepsilon $ & $\varepsilon $ & $\varepsilon $ \\ 
$\alpha $ & $\rho $ & $\varepsilon $ & $\rho $ & $\varepsilon $ & $\rho $ & $%
\varepsilon $ & $\rho $ & $\varepsilon $ \\ 
$\beta $ & $\upsilon $ & $\upsilon $ & $\varepsilon $ & $\varepsilon $ & $%
\upsilon $ & $\upsilon $ & $\varepsilon $ & $\varepsilon $ \\ 
$\gamma $ & $\tau $ & $\upsilon $ & $\rho $ & $\varepsilon $ & $\tau $ & $%
\upsilon $ & $\rho $ & $\varepsilon $ \\ 
$\tau $ & $\gamma $ & $\gamma $ & $\gamma $ & $\gamma $ & $\varepsilon $ & $%
\varepsilon $ & $\varepsilon $ & $\varepsilon $ \\ 
$\upsilon $ & $\beta $ & $\gamma $ & $\beta $ & $\gamma $ & $\rho $ & $%
\varepsilon $ & $\rho $ & $\varepsilon $ \\ 
$\rho $ & $\alpha $ & $\alpha $ & $\gamma $ & $\gamma $ & $\upsilon $ & $%
\upsilon $ & $\varepsilon $ & $\varepsilon $ \\ 
$\varepsilon $ & $0$ & $\alpha $ & $\beta $ & $\gamma $ & $\tau $ & $%
\upsilon $ & $\rho $ & $\varepsilon $%
\end{tabular}%
\text{ }%
\begin{tabular}{l|llllllll}
$\oplus _{21}^{8}$ & $0$ & $\alpha $ & $\beta $ & $\gamma $ & $\tau $ & $%
\upsilon $ & $\rho $ & $\varepsilon $ \\ \hline
$0$ & $0$ & $\alpha $ & $\beta $ & $\gamma $ & $\tau $ & $\upsilon $ & $\rho 
$ & $\varepsilon $ \\ 
$\alpha $ & $\alpha $ & $\alpha $ & $\gamma $ & $\gamma $ & $\upsilon $ & $%
\upsilon $ & $\varepsilon $ & $\varepsilon $ \\ 
$\beta $ & $\beta $ & $\gamma $ & $\beta $ & $\gamma $ & $\rho $ & $%
\varepsilon $ & $\rho $ & $\varepsilon $ \\ 
$\gamma $ & $\gamma $ & $\gamma $ & $\gamma $ & $\gamma $ & $\varepsilon $ & 
$\varepsilon $ & $\varepsilon $ & $\varepsilon $ \\ 
$\tau $ & $\tau $ & $\upsilon $ & $\rho $ & $\varepsilon $ & $\tau $ & $%
\upsilon $ & $\rho $ & $\varepsilon $ \\ 
$\upsilon $ & $\upsilon $ & $\upsilon $ & $\varepsilon $ & $\varepsilon $ & $%
\upsilon $ & $\upsilon $ & $\varepsilon $ & $\varepsilon $ \\ 
$\rho $ & $\rho $ & $\varepsilon $ & $\rho $ & $\varepsilon $ & $\rho $ & $%
\varepsilon $ & $\rho $ & $\varepsilon $ \\ 
$\varepsilon $ & $\varepsilon $ & $\varepsilon $ & $\varepsilon $ & $%
\varepsilon $ & $\varepsilon $ & $\varepsilon $ & $\varepsilon $ & $%
\varepsilon $%
\end{tabular}
\tag{3.2.}
\end{equation}

In this case, we remark that in this MV-algebra structure all elements are
idempotent, therefore it is a Boolean algebra. For $\alpha $ and $\rho $ we
remark that $\alpha \vee \beta \vee \tau =\varepsilon $, $\alpha \wedge
\beta \wedge \tau =0$ and $\upsilon \vee \rho \vee \gamma =\varepsilon $, $%
\upsilon \wedge \rho \wedge \gamma =0$. Therefore, as MV-algebra, $W$ $%
\simeq (-\infty ,\alpha ]\times (-\infty ,\beta ]\times (-\infty ,\tau ]$ or 
$W$ $\simeq (-\infty ,\upsilon ]\times (-\infty ,\rho ]\times (-\infty
,\gamma ]$, where $(-\infty ,\alpha ]=\{0,\alpha \}$, $(-\infty ,\beta
]=\{0,\beta \}$, $(-\infty ,\gamma ]=\{0,\gamma \}$, $(-\infty ,\tau
]=\{0,\tau \}$, $(-\infty ,\upsilon ]=\{O,\upsilon \}$, $(-\infty ,\rho
]=\{0,\rho \}$.

From here, we get that there are only three non-isomorphic MV-algebras of
order $8$. Thus, it results that in an MV-algebra of order $8$ we can have
only $0$, $2$ or $6$ proper idempotents.\medskip 
\begin{equation*}
\end{equation*}

\bigskip \textbf{4. Binary block codes associated to a Boolean algebra}

\begin{equation*}
\end{equation*}

In this section, starting from a Boolean algebra of order $2$, we will give
an algorithm to built a Boolean algebra of order $2^{k+1},k\geq 1$. We
denote such an algebra by $\mathcal{B}_{2^{k+1}}$.

In [FV; 20], to an MV-algebra and to a Wajsberg algebra were associated
binary block codes and, in some circumstances, was proved that the converse
is also true. Using some of these ideas, to algebra $\mathcal{B}_{2^{k+1}}~$%
we will associate a binary block code and we will prove that the converse of
this statement is also true, namely to such a binary block code a Boolean
algebra $\mathcal{B}_{2^{k+1}}~$can be associated.\medskip 

\textbf{Definition 4.1.} Two Boolean algebras $\left( \mathcal{B},\vee
\wedge ,\rceil ,0,1\right) $ and $\left( \mathcal{B}^{\prime },\curlyvee
,\curlywedge ,\widetilde{},\mathbf{0},\mathbf{1}\right) $ are said to be 
\textit{isomorphic} if and only if there is a bijective function $f:\mathcal{%
B}\rightarrow \mathcal{B}^{\prime }$ satisfying the following conditions:

i) $f\left( x\vee y\right) =f\left( x\right) \vee f\left( y\right) ,$for all 
$x,y\in \mathcal{B};$

ii) $f\left( x\wedge y\right) =f\left( x\right) \wedge f\left( y\right) ,$
for all $x,y\in \mathcal{B};$

iii) $f\left( \rceil x\right) =\widetilde{f\left( x\right) }$, for all $x\in 
\mathcal{B};$

iv) $f\left( 0\right) =\mathbf{0};$

v) $f\left( 1\right) =\mathbf{1}.\medskip $

Let $\mathcal{B}_{2}$ be the boolean algebra given in the following table

\begin{equation}
\begin{tabular}{|l|l|l|}
\hline
$\oplus _{11}^{2}$ & $\beta $ & $\varepsilon $ \\ \hline
$\beta $ & $\beta $ & $\varepsilon $ \\ \hline
$\varepsilon $ & $\varepsilon $ & $\varepsilon $ \\ \hline
\end{tabular}%
.  \tag{4.1.}
\end{equation}

To Boolean algebra $\mathcal{B}_{2}$ we will attach the map 
\begin{equation}
\varphi _{2}:\mathcal{B}_{2}\times \mathcal{B}_{2}\rightarrow \mathcal{B}%
_{2},\varphi _{2}\left( x,y\right) =x\oplus _{11}^{2}y.  \tag{4.2.}
\end{equation}

In the following, we will use the same notation $\mathcal{B}_{2}$ for the
table 
\begin{tabular}{|l|l|}
\hline
$\beta $ & $\varepsilon $ \\ \hline
$\varepsilon $ & $\varepsilon $ \\ \hline
\end{tabular}%
.

We consider $\mathcal{C}_{2}$ a Boolean algebra of order $2$ isomorphic to $%
\mathcal{B}_{2},$ which has the following multiplication table%
\begin{equation}
\begin{tabular}{|l|l|l|}
\hline
$\oplus _{11}^{\prime 2}$ & $0$ & $\alpha $ \\ \hline
$0$ & $0$ & $\alpha $ \\ \hline
$\alpha $ & $\alpha $ & $\alpha $ \\ \hline
\end{tabular}%
.  \tag{4.3.}
\end{equation}

\bigskip

For table 
\begin{tabular}{|l|l|}
\hline
$0$ & $\alpha $ \\ \hline
$\alpha $ & $\alpha $ \\ \hline
\end{tabular}%
, we use the same notation $\mathcal{C}_{2}$. Let $f_{2}:\mathcal{C}%
_{2}\rightarrow \mathcal{B}_{2}$ be an isomorphism of Boolean algebras.

To Boolean algebra $\mathcal{C}_{2}$ we will attach the map 
\begin{equation}
\theta _{2}:\mathcal{C}_{2}\times \mathcal{C}_{2}\rightarrow \mathcal{C}%
_{2},\theta _{2}\left( x,y\right) =x\oplus _{11}^{\prime 2}y.  \tag{4.4.}
\end{equation}%
Therefore, the multiplication table of $\mathcal{B}_{4}$, given in $\left(
3.1\right) $, can be written under the form 
\begin{equation*}
\mathcal{B}_{4}=%
\begin{tabular}{|l|l|}
\hline
$\mathcal{C}_{2}$ & $\mathcal{B}_{2}$ \\ \hline
$\mathcal{B}_{2}$ & $\mathcal{B}_{2}$ \\ \hline
\end{tabular}%
.
\end{equation*}%
We \ remark that $\mathcal{B}_{4}=\mathcal{C}_{2}\cup \mathcal{B}_{2}$ and $%
\mathcal{C}_{2}\cap \mathcal{B}_{2}=\varnothing $. To Boolean algebra $%
\mathcal{B}_{4}$ we will attach the map 
\begin{equation}
\varphi _{4}:\mathcal{B}_{4}\times \mathcal{B}_{4}\rightarrow \mathcal{B}%
_{4},\varphi _{4}\left( x,y\right) \text{=}\left\{ 
\begin{array}{c}
\theta _{2}\left( x,y\right) ,\text{ for }x,y\in \mathcal{C}_{2} \\ 
\varphi _{2}\left( f_{2}(x),y\right) ,\text{ for }x\in \mathcal{C}_{2},y\in 
\mathcal{B}_{2} \\ 
\varphi _{2}\left( x,f_{2}(y)\right) ,\text{ for }x\in \mathcal{B}_{2},y\in 
\mathcal{C}_{2} \\ 
\varphi _{2}\left( x,y\right) ,\text{ for }x\in \mathcal{B}_{2},y\in 
\mathcal{B}_{2}%
\end{array}%
\right. .  \tag{4.5.}
\end{equation}%
From the above, it is easy to see that $\varphi _{4}\left( x,y\right)
=x\oplus _{11}^{4}y$.

Continuing the algorithm, with the above notations, we remark that the
multiplication table of the Boolean algebra $\mathcal{B}_{8}$ can be written
under the form%
\begin{equation*}
\mathcal{B}_{8}=%
\begin{tabular}{|l|l|}
\hline
$\mathcal{C}_{4}$ & $\mathcal{B}_{4}$ \\ \hline
$\mathcal{B}_{4}$ & $\mathcal{B}_{4}$ \\ \hline
\end{tabular}%
,
\end{equation*}%
where $(\mathcal{C}_{4},\oplus _{11}^{\prime 4})$ is a Boolean algebra
isomorphic to $\mathcal{B}_{4}$. Let $f_{4}:\mathcal{C}_{4}\rightarrow 
\mathcal{B}_{4}$ be an isomorphism of Boolean algebras. To Boolean algebra $%
\mathcal{C}_{4}$ we will attach the map 
\begin{equation*}
\theta _{4}:\mathcal{C}_{4}\times \mathcal{C}_{4}\rightarrow \mathcal{C}%
_{4},\theta _{4}\left( x,y\right) =x\oplus _{11}^{\prime 4}y.
\end{equation*}%
We \ remark that $\mathcal{B}_{8}=\mathcal{C}_{4}\cup \mathcal{B}_{4}$ and $%
\mathcal{C}_{4}\cap \mathcal{B}_{4}=\varnothing $. To Boolean algebra $%
\mathcal{B}_{4}$ we will attach the map 
\begin{equation*}
\varphi _{8}:\mathcal{B}_{8}\times \mathcal{B}_{8}\rightarrow \mathcal{B}%
_{8},\varphi _{8}\left( x,y\right) \text{=}\left\{ 
\begin{array}{c}
\theta _{4}\left( x,y\right) ,\text{ for }x,y\in \mathcal{C}_{4} \\ 
\varphi _{4}\left( f_{4}(x),y\right) ,\text{ for }x\in \mathcal{C}_{4},y\in 
\mathcal{B}_{4} \\ 
\varphi _{4}\left( x,f_{4}(y)\right) ,\text{ for }x\in \mathcal{B}_{4},y\in 
\mathcal{C}_{4} \\ 
\varphi _{4}\left( x,y\right) ,\text{ for }x,y\in \mathcal{B}_{4}%
\end{array}%
\right. .
\end{equation*}%
From here, it is easy to see that $\varphi _{8}\left( x,y\right) =x\oplus
_{21}^{8}y$.

Therefore, following the above algorithm, the multiplication table of the
Boolean algebra $\mathcal{B}_{2^{k}}$ can be written under the form%
\begin{equation}
\mathcal{B}_{2^{k}}=%
\begin{tabular}{|l|l|}
\hline
$\mathcal{C}_{2^{k-1}}$ & $\mathcal{B}_{2^{k-1}}$ \\ \hline
$\mathcal{B}_{2^{k-1}}$ & $\mathcal{B}_{2^{k-1}}$ \\ \hline
\end{tabular}%
,  \tag{4.6.}
\end{equation}%
where $\mathcal{C}_{2^{k-1}}$ is a Boolean algebra isomorphic to $\mathcal{B}%
_{2^{k-1}}$.\medskip

\textbf{Algorithm 1}

Assuming that we built the Boolean algebra $(\mathcal{B}_{2^{k}},\oplus
_{21}^{2^{k}})$ with $\varphi _{2^{k}}:\mathcal{B}_{2^{k}}\times \mathcal{B}%
_{2^{k}}\rightarrow \mathcal{B}_{2^{k}},\varphi _{2^{k}}\left( x,y\right)
=x\oplus _{21}^{2^{k}}y$ and considering $(\mathcal{C}_{2^{k}},\oplus
_{11}^{\prime 2^{k}})$ a Boolean algebra isomorphic to $\mathcal{B}_{2^{k}}$%
, let $f_{2^{k}}:\mathcal{C}_{2^{k}}\rightarrow \mathcal{B}_{_{2^{k}}}$ be
an isomorphism of Boolean algebras.

To Boolean algebra $\mathcal{C}_{2^{k}}$, we will attach the map 
\begin{equation}
\theta _{2^{k}}:\mathcal{C}_{2^{k}}\times \mathcal{C}_{2^{k}}\rightarrow 
\mathcal{C}_{2^{k}},\theta _{2^{k}}\left( x,y\right) =x\oplus _{11}^{\prime
2^{k}}y.  \tag{4.7.}
\end{equation}

We consider the set $\mathcal{B}_{2^{k+1}}=\mathcal{C}_{2^{k}}\cup \mathcal{B%
}_{2^{k}}$ with $\mathcal{C}_{2^{k}}\cap \mathcal{B}_{2^{k}}=\varnothing $.
To Boolean algebra $\mathcal{B}_{2^{k+1}}$ we will attach the map 
\begin{equation}
\varphi _{2^{k+1}}:\mathcal{B}_{2^{k+1}}\times \mathcal{B}%
_{2^{k+1}}\rightarrow \mathcal{B}_{2^{k+1}},\varphi _{2^{k+1}}\left(
x,y\right) \text{=}\left\{ 
\begin{array}{c}
\theta _{2^{k}}\left( x,y\right) ,\text{ for }x,y\in \mathcal{C}_{2^{k}} \\ 
\varphi _{2^{k}}\left( f_{2^{k}}(x),y\right) ,\text{ for }x\in \mathcal{C}%
_{2^{k}},y\in \mathcal{B}_{2^{k}} \\ 
\varphi _{2^{k}}\left( x,f_{2^{k}}(y)\right) ,\text{ for }x\in \mathcal{B}%
_{2^{k}},y\in \mathcal{C}_{2^{k}} \\ 
\varphi _{2^{k}}\left( x,y\right) ,\text{ for }x\in \mathcal{B}_{2^{k}},y\in 
\mathcal{B}_{2^{k}}%
\end{array}%
\right. .  \tag{4.8.}
\end{equation}%
Defining $x\oplus _{21}^{2^{k+1}}y=\varphi _{2^{k+1}}\left( x,y\right) $, it
results that $(\mathcal{B}_{2^{k+1}},\oplus _{21}^{2^{k+1}})$ is a Boolean
algebra of order $2^{k+1}$, as can easily be checked.

Let $\left( X,\oplus ,\lceil ,0\right) $ be a finite MV-algebra of order $n$%
, with $X=\{0=\alpha _{0},\alpha _{1},\alpha _{2},...,\alpha
_{n-2},\varepsilon =\alpha _{n-1}\}$. In [FV; 20], to an MV-algebra and to
associated Wajsberg algebra were associated binary block codes. Summarizing
these methods, we give the following definition. Let $C_{n}=%
\{w_{0},w_{1},...,w_{n-2},w_{\varepsilon }\}$ be a binary block code, with
codewords of length $n$.\medskip 

\textbf{Definition 4.2.} 1) The block code $C$ is attached to MV-algebra $X,$
if for a codeword $w_{j}\in X$, $w_{j}=i_{0}i_{1}...i_{n-2}i_{\varepsilon
},i_{0},i_{1},...i_{n-2},i_{\varepsilon }\in \{0,1\},j\in
\{0,1,2,...,n-2,\varepsilon \},$ we have that $i_{s}=1$ if $\alpha
_{j}\oplus \alpha _{s}=\varepsilon $ and $i_{s}=0,$otherwise, $s\in
\{0,1,2,...,n-2,\varepsilon \}$.

2) A matrix attached to the code $C$, is a quadratic matrix $M_{C}=\left(
m_{i,j}\right) _{i,j\in \{1,2,...,n\}}\in \mathcal{M}_{n}(\{0,1\})$ such
that its rows are formed by the codewords of $C.\medskip $

\textbf{Remark 4.3.} Since a Boolean algebra is an MV-algebra, we have that:

- the code $C_{2}$ attached the the algebra $\mathcal{B}_{2}$ is $%
C_{2}=\{01,11\}=\{w_{0},w_{1}\}$ and the attached matrix is 
\begin{equation*}
M_{C_{2}}=\left( 
\begin{array}{cc}
0 & 1 \\ 
1 & 1%
\end{array}%
\right) .
\end{equation*}

- the code $C_{4}$ attached to the algebra $\mathcal{B}_{4}$ is $%
C_{4}=\{0001,0011,0101,1111\}=\{w_{0},w_{1},w_{2},w_{3}\}$ and the attached
matrix is 
\begin{equation*}
M_{C_{4}}=\left( 
\begin{array}{cccc}
0 & 0 & 0 & 1 \\ 
0 & 0 & 1 & 1 \\ 
0 & 1 & 0 & 1 \\ 
1 & 1 & 1 & 1%
\end{array}%
\right) .
\end{equation*}

- the code attached to the algebra $\mathcal{B}_{8}$ is $C_{8}=%
\{00000001,00000011,00000111,$\newline
$00001111,00010001,00110011,01110111,11111111\}$ and the attached matrix is 
\begin{equation*}
M_{C_{8}}=\left( 
\begin{array}{cccccccc}
0 & 0 & 0 & 0 & 0 & 0 & 0 & 1 \\ 
0 & 0 & 0 & 0 & 0 & 0 & 1 & 1 \\ 
0 & 0 & 0 & 0 & 0 & 1 & 0 & 1 \\ 
0 & 0 & 0 & 0 & 1 & 1 & 1 & 1 \\ 
0 & 0 & 0 & 1 & 0 & 0 & 0 & 1 \\ 
0 & 0 & 1 & 1 & 0 & 0 & 1 & 1 \\ 
0 & 1 & 0 & 1 & 0 & 1 & 0 & 1 \\ 
1 & 1 & 1 & 1 & 1 & 1 & 1 & 1%
\end{array}%
\right) .
\end{equation*}

\bigskip If we denote with $\mathbf{0}_{n}$ the zero matrix with $n$
elements, we remark that 
\begin{equation*}
M_{C_{4}}=\left( 
\begin{array}{cc}
\mathbf{0}_{2} & M_{C_{2}} \\ 
M_{C_{2}} & M_{C_{2}}%
\end{array}%
\right) ,M_{C_{8}}=\left( 
\begin{array}{cc}
\mathbf{0}_{4} & M_{C_{4}} \\ 
M_{C_{4}} & M_{C_{4}}%
\end{array}%
\right) .
\end{equation*}%
Therefore, we have 
\begin{equation}
M_{C_{2^{k+1}}}=\left( 
\begin{array}{cc}
\mathbf{0}_{_{2^{k}}} & M_{C_{_{2^{k}}}} \\ 
M_{C_{_{2^{k}}}} & M_{C_{_{2^{k}}}}%
\end{array}%
\right) ,  \tag{4.9.}
\end{equation}%
with $C_{2^{k+1}}$ the attached binary block code having as codewords the
rows of the matrix $M_{C_{2^{k+1}}}.$

In the following, we provide a method to attach to a binary block code $%
C_{2^{k+1}}$ a Boolean algebra.

First of all, we consider the binary block code $C_{2}=\{w_{0}\preceq
w_{1}\},$ with $\preceq $ the lexicographic order. We define the following
multiplication:%
\begin{equation}
w_{0}\ast _{2}w_{1}=w_{1}\ast _{2}w_{0}=w_{1}\ast _{2}w_{1}=w_{1},\text{ }%
w_{0}\ast _{2}w_{0}=w_{0}.  \tag{4.10.}
\end{equation}%
It results that $\left( C_{2},\ast \right) $ is a Boolean algebra,
isomorphic to $\mathcal{B}_{2}$.

If we consider $C_{4}=\{w_{0}\preceq w_{1}\preceq w_{2}\preceq w_{3}\}$ with
the lexicographic order $\preceq $, let $C_{4}^{\prime }=\{w_{0}\preceq
w_{1}\}$ and $C_{4}^{\prime \prime }=\{w_{2}\preceq w_{3}\}$ be two disjoint
subsets of $C_{4}$. We have that $(C_{4}^{\prime },\ast _{2})$ is a Boolean
algebra of order $2$. On $C_{4}^{\prime \prime }$ we define the
multiplication

\begin{equation*}
w_{2}\ast _{2}^{\prime }w_{3}=w_{3}\ast _{2}^{\prime }w_{2}=w_{3}\ast
_{2}^{\prime }w_{3}=w_{3},\text{ }w_{2}\ast _{2}^{\prime }w_{2}=w_{2}.
\end{equation*}%
It results that $(C_{4}^{\prime \prime },\ast _{2}^{\prime })$ is a Boolean
algebra isomorphic to $(C_{4}^{\prime },\ast _{2})$. On $C_{4}$ we define
the following multiplication 
\begin{equation*}
w_{i}\ast _{4}w_{j}\text{=}\left\{ 
\begin{array}{c}
w_{i}\ast _{2}w_{j},\text{ for }w_{i},w_{j}\in C_{4}^{\prime }, \\ 
w_{i}\ast _{2}^{\prime }w_{j},\text{ for }w_{i}\in C_{4}^{\prime },w_{j}\in
C_{4}^{\prime \prime }, \\ 
w_{i}\ast _{4}^{\prime }w_{j},\text{ for }w_{i}\in C_{4}^{\prime \prime
},w_{j}\in C_{4}^{\prime }, \\ 
w_{i}\ast _{4}^{\prime }w_{j},\text{ for }w_{i},w_{j}\in C_{4}^{\prime
\prime }.%
\end{array}%
\right. .
\end{equation*}%
From here, we have that $(C_{4},\ast _{4})$ is a Boolean algebra isomorphic
to $\mathcal{B}_{4}$.\medskip

\textbf{Algorithm 2}

Assuming that we have defined the Boolean algebra $\left( C_{2^{k}},\ast
_{2^{k}}\right) $ isomorphic to $\mathcal{B}_{2^{k}}$, let $%
C_{2^{k+1}}=\{w_{0}\preceq w_{1}\preceq ,...\preceq w_{2^{k}-1}\preceq
...\preceq w_{2^{k+1}-1}\}$ be the binary block code defined by the matrix $%
M_{C_{2^{k+1}}}$, given by the relation $\left( 4.9\right) $, with the
codewords lexicographically ordered. We consider the sets $%
C_{2^{k+1}}^{\prime }=\{w_{0}\preceq w_{1}\preceq ,...\preceq w_{2^{k}-1}\}$
and $C_{2^{k+1}}^{\prime \prime }=\{w_{2^{k}}\preceq ...\preceq
w_{2^{k+1}-1}\}$. We have that $(C_{2^{k+1}}^{\prime },\ast _{2^{k}})$ is
isomorphic to $\mathcal{B}_{2^{k}}$ and on $C_{2^{k+1}}^{\prime \prime }$ we
define a multiplication $\ast _{2^{k}}^{\prime }$ such that \ also $%
(C_{2^{k+1}}^{\prime \prime },\ast _{2^{k}}^{\prime })$ to be isomorphic to $%
\mathcal{B}_{2^{k}}$.

On $C_{2^{k+1}}$ we define the following multiplication 
\begin{equation}
w_{i}\ast _{2^{k+1}}w_{j}\text{=}\left\{ 
\begin{array}{c}
w_{i}\ast _{2^{k+1}}w_{j},\text{ for }w_{i},w_{j}\in C_{2^{k+1}}^{\prime },
\\ 
w_{i}\ast _{2^{k+1}}^{\prime }w_{j},\text{ for }w_{i}\in C_{2^{k+1}}^{\prime
},w_{j}\in C_{2^{k+1}}^{\prime \prime }, \\ 
w_{i}\ast _{2^{k+1}}^{\prime }w_{j},\text{ for }w_{i}\in C_{2^{k+1}}^{\prime
\prime },w_{j}\in C_{2^{k+1}}^{\prime }, \\ 
w_{i}\ast _{2^{k+1}}^{\prime }w_{j},\text{ for }w_{i},w_{j}\in
C_{2^{k+1}}^{\prime \prime }.%
\end{array}%
\right. .  \tag{4.11.}
\end{equation}%
It results that $(C_{2^{k+1}},\ast _{2^{k+1}})$ is a Boolean algebra
isomorphic to $\mathcal{B}_{2^{k+1}}$.\medskip

From the above results, we proved the following Theorem.\medskip 

\textbf{Theorem 4.4. }

\textit{1) To each Boolean algebra of order} $2^{k+1}$, $\mathcal{B}%
_{2^{k+1}}$, \textit{we can associate a binary block code} $C_{2^{k+1}}$, 
\textit{with associated matrix given by the relation} $\left( 4.9\right) $.

2) \textit{On binary block code} $C_{2^{k+1}}$ \textit{we can define a
multiplication} $\ast _{2^{k+1}}$\textit{such that} $(C_{2^{k+1}},\ast
_{2^{k+1}})$ \textit{is a Boolean algebra isomorphic to} $\mathcal{B}%
_{2^{k+1}}$. $_{{}}\medskip $ 
\begin{equation*}
\end{equation*}

\textbf{Conclusions.} In this paper, we gave some properties of MV-algebras.
We defined a Fibonacci sequence in such an algebra and we proved that a
Fibonacci stationary sequence gives us an idempotent element. Taking into
account of the representation of a finite MV-algebra, by using Boolean
elements of this algebra, we proved that such a sequence is always
stationary. This result is interesting comparing with the behavior of such a
sequence on the group $\left( \mathbb{Z}_{n},+\right) $, the group of
integers modulo $n$, \ where the Fibonacci sequences are periodic, with
period given by Pisano period.

In Section 3, some examples of MV-algebras and the number of their
idempotents are given. In Section 4, was provided an algorithm to built a
Boolean algebra of order $2^{k+1}$, starting from a Boolean algebra of order
two. Moreover, as an application in Coding Theory, to a Boolean algebra was
attached a binary block code and was proved that, under some conditions, the
converse is also true.\medskip\ 
\begin{equation*}
\end{equation*}

\textbf{Acknowledgements}. This work continues some of the results presented
at IECMSA-2019 conference. The author thanks organizers of IECMSA-2019,
especially Professor Murat Tosun and Professor Soley Ersoy for their support.

\bigskip

\begin{equation*}
\end{equation*}

\textbf{References}%
\begin{equation*}
\end{equation*}

[BV; 10] Belohlavek, R., Vychodil, V., \textit{Residuated Lattices of Size} $%
\leq 12$, Order, 27(2010), 147-161.

[CHA; 58] Chang, C.C.,\textit{\ Algebraic analysis of many-valued logic},
Trans. Amer. Math. Soc. 88(1958), 467-490.

[COM; 00] Cignoli, R. L. O, Ottaviano, I. M. L. D, Mundici, D., \textit{%
Algebraic foundations of many-valued reasoning}, Trends in Logic, Studia
Logica Library, Dordrecht, Kluwer Academic Publishers, 7(2000).

[FHSV; 20] Flaut, C., Hoskova-Mayerova, S., Saeid, A.B., Vasile, R., \textit{%
Wajsberg algebras of order n~}$(n\leq 9)$, accepted in Neural Computing and
Applications, DOI: 10.1007/s00521-019-04676-x.

[FV; 20] Flaut, C.,Vasile, R., \textit{Wajsberg algebras arising from binary
block codes}, accepted in Soft Computing, DOI: 10.1007/s00500-019-04653-5.

[FRT; 84] Font, J., M., Rodriguez, A., J., Torrens, A., \textit{Wajsberg
Algebras}, Stochastica, 8(1)(1984), 5-30.

[HKN; 12] Han, J.S., Kim, H.S., Neggers, J., \textit{Fibonacci sequences in
groupoids}, Advances in Difference Equations 2012, 2012:19.

[HR; 99] H\"{o}hle, U., Rodabaugh, S., E., \textit{Mathematics of Fuzzy
Sets: Logic, Topology and Measure Theory}, Springer Science and Business
Media, LLC, 1999.

[KNS; 13] Kim, H.S., Neggers, J., So, K.S., \textit{Generalized Fibonacci
sequences in groupoids}, Advances in Difference Equations 2013, 2013:26.

[Mu; 07] Mundici, D., \textit{MV-algebras-a short tutorial}, Department of
Mathematics \textit{Ulisse Dini}, University of Florence, 2007.

[Re; 13] M. Renault, \textit{The Period, Rank, and Order of the (a,
b)-Fibonacci Sequence Mod m}, Mathematics Magazine, 86(5)(2013), 372-380,

https://doi.org/10.4169/math.mag.86.5.372.

[WDH; 17] Wang, J., T., Davvaz, B., He, P., F., \textit{On derivations of
MV-algebras}, https://arxiv.org/pdf/1709.04814.pdf.

\begin{equation*}
\end{equation*}

Cristina Flaut

{\small Faculty of Mathematics and Computer Science, Ovidius University,}

{\small Bd. Mamaia 124, 900527, CONSTANTA, ROMANIA}

{\small \ http://www.univ-ovidius.ro/math/}

{\small e-mail: cflaut@univ-ovidius.ro; cristina\_flaut@yahoo.com}

\end{document}